\documentclass[review]{elsarticle}

\usepackage{lineno,hyperref}
\modulolinenumbers[5]

\journal{Journal of \LaTeX\ Templates}

\usepackage{amsfonts,amssymb,amsbsy}
\usepackage{latexsym}
\usepackage{amsmath}
\usepackage{xcolor}

\newtheorem{theorem}{Theorem}
\newtheorem{lemma}{Lemma}
\newtheorem{corollary}{Corollary}
\newtheorem{example}{Example}
\newtheorem{remark}{Remark}
\newtheorem{definition}{Definition}

\newcommand{\R}{\mathbb{R}}

\newcommand{\Z}{\mathbb{Z}}
\begin{document}
	
	\begin{frontmatter}
		
		\title{Absolute exponential stability criteria of delay time-varying  systems with sector-bounded nonlinearity: a comparison approach }
		
		\author[mymainaddress]{Nguyen Khoa Son\corref{mycorrespondingauthor}}
		\cortext[mycorrespondingauthor]{Corresponding author}
		\address[mymainaddress]{Institute of Mathematics, Vietnam Academy of Science and Technology, 18 Hoang Quoc Viet Rd., Hanoi, Vietnamm}\ead{nkson@vast.vn}
		
		\author[mysecondaryaddress]{Nguyen Thi Hong}
		
		\address[mysecondaryaddress]{Institute of Mathematics, Vietnam Academy of Science and Technology, 18 Hoang Quoc Viet Rd., Hanoi, Vietnamm}
		\ead{nthong@math.ac.vn}
		
		
		\begin{abstract}
		Absolute exponential stability problem of  delay time-varying systems (DTVS) with sector-bounded nonlinearity is presented in this paper. By using  the comparison principle and properties of positive systems we derive several novel  criteria of absolute exponential stability, for both continuous-time and discrete-time nonlinear DTVS. When applied to the time-invariant case,  the obtained stability criteria are shown to cover and extend some previously known results, including, in particular, the result due to S.K. Persidskii in Ukrainian Mathematical Journal,  vol. 57(2005). The theoretical results are illustrated by examples that can not be treated by the existing ones. 
\end{abstract}

	\end{frontmatter}

\section{Introduction and Preliminaries}
\label{sec:introduction}
The absolute stability problem, first formulated in \cite{lurie1944}, is one of the main problems in the systems and control theory. Roughly speaking, a dynamic nonlinear system containing nonlinearities in its mathematical description  is said to be absolute stable if its equilibrium is asymptotically stable for any nonlinearity in a given nonlinearities class. In this nonlinear framework, the most widely used approach for characterizing and checking the absolute stability  is the Lyapunov function method. The reader is referred to the survey article \cite{liberzon_automat2008} and the monograph \cite{liao_book} for the study on absolute stability of the so-called Luri'e control systems and the monograph \cite{KaszBhaya} which is dedicated to absolute stability analysis of several classes of nonlinear systems of practical interest, including the so-called Persidskii-type systems, whose stability can be characterized by diagonal Lyapunov functions. 

As is well-known, the most simple Persidskii-type system can be represented by the differential equation
\begin{equation}\label{persid1}
	\dot x= Af(x), \ f(0)=0 	
\end{equation}
where $A:=(a_{ij})$ is $n\times n$-matrix and  $f$ is supposed to be  a  continuous diagonal function  $f(x):=(f_1(x_1),\ldots, f_n(x_n))^{\top}$ belonging to the {\it infinite sector}
\begin{equation}\label{infsector}
K(0,\infty):= \{f: \; x_if_i(x_i) > 0, \forall x_i\not=0,\ i=1,\ldots, n\}.	
\end{equation}
This class of models was first introduced for stability analysis in \cite{barbashin},  where a linear combination of the integrals of the nonlinearities was used as a Lyapunov
function. Next, that result was improved and extended by Persidskii \cite{persidskii69}   and  in many subsequent papers (see, e.g. \cite{K_B-SIAM1993, Oliveira2002, persidskii2005,   Efimov2021}), where different types of Lyapunov functions  have been proposed  for checking absolute stability for continuous-time systems as well as their discrete-time counterparts, under different classes of sector nonlinearities. The study of absolute stability of Persidskii-type systems is an
important tool in many application fields such as automatic control, Lotka-Volterra ecosystems,  Hopfield neural networks,  and decentralized
power-frequency control of power systems, among others (see, e.g. \cite{K_B-SIAM1993}). Recently, the absolute stability problems have been investigated  also for the classes of switched nonlinear systems, including those with time-delay in the state variables, see, e.g. \cite{ sun_wang2013, alex_mason2014},  \cite{zhang_zhao2016} \cite{alex2021} and the references therein.  It is important to note that in most of the aforementioned works only the classes of time-invariant nonlinear systems have been considered where some criteria of absolute asymptotical stability have been derived. So far little  attention has been devoted to time-varying systems and  absolute exponential stability analysis.

In this paper we will study the absolute exponential stability problem for the Persidskii class of time-varying delay nonlinear systems of the form
\begin{equation} \label{TVS}
	\dot x=A(t) f(x(t)) + B(t)f(x(t-h)), \ t\geq 0	
	\end{equation}
where $A(\cdot), B(\cdot)$ are $n\times n$-matrix continuous functions and the diagonal function $f$ belongs to  the {\it bounded sector} $K[\delta,\beta]$ defined as 
\begin{equation}\label{sector2}
K[\delta,\beta]\!:=\! \{f: \delta_ix_i^2\leq x_if_i(x_i) \leq \beta_i x_i^2,\forall x_i\not=0, \ i=1,\ldots, n	\}
\end{equation}
with $0<\delta_i\leq \beta_i, i\in \underline n$ being given numbers. Our primary purpose  is to derive some verifiable  criteria of absolute exponential stability  for this class of time-varying nonlinear systems. Similar results will be established also for time-varying difference  systems.  What is more, it is remarkable that when applied to the particular cases of time-invariant systems, the obtained results yield novel criteria of absolute exponential stability which improve the existing ones. Differently from the previous works which are mainly based on the Lyapunov functions (or Lyapunov-Krasovskii functions) method, we will employ the comparison principle in deriving the main results. It is worthy to mention that the comparison approach was proved to be a quite effective tool in stability analysis, particularly for time-varying and switched systems. The readers are referred to \cite{zhang2012, Liu2018,Tian_Sun,son_ngoc_2022} and the literature given therein on this topic. Basically, this approach is based on comparing the original systems with a positive comparison system of the same dimension and using the stability characterization of the latter one (which is supported by the powerful Perron-Frobenius theory, see e.g. \cite{Chelaboina2004,  Ngoc_IEEE, blanc_valcher}) to conclude on stability of the original system. 

\noindent {\it Notation.}
Throughout the paper, $\R, \R_+,\Z, \Z_+$  will stand for the sets of real numbers, nonnegative real numbers,   natural numbers and nonnegative natural numbers,  respectively. For $t_1,t_2 \in \Z$ with $ t_1<t_2, \Z_{[t_1,t_2]}$ denotes the interval in $\Z: \Z_{[t_1,t_2]}= \{t_1,t_1+1,\ldots, t_2-1,t_2\}$.  For an integer $m, \underline m$ denotes the set of numbers $ \{1,2,\ldots, m \}$. Inequalities between vectors and matrices are understood componentwise: for  vectors $x=(x_i), y=(y_i) \in \R^n$  we write $x\geq y$ and $x\gg y$ iff $x_i \geq y_i$ and $x_i> y_i$, for all $i\in \underline n,$ respectively. Denote $|x|=(|x_i|)$ and $x^{\top}$ is the transpose of $x$. Similar notation is adopted for $(n\times n)$-matrices. Without loss of generality, the norm of vectors $x\in \R^n$ is assumed to be the $1$-norm: $\|x\|=\sum_{i=1}^n|x_i|$. A matrix $A\in \R^{n\times n}$ is said to be {\it  Hurwitz stable} if  $\text{Re}\lambda <0$ and {\it Schur stable} if $|\lambda|<1$,  for any root  $\lambda $ of the characteristic polynomial, i.e. $  \det(\lambda I-A)=0$.  Matrix $A=(a_{ij})\in \R^{n\times n}$ is said to be a Metzler matrix if $a_{ij}\geq 0$ for all $i\not=j$. The following stability characterization of Metzler and nonnegative matrices is useful in stability analysis of positive systems,  see e.g. \cite{Berman}.
\begin{lemma}\label{metzmatrix} Assume that $A=(a_{ij})\in \R^{n\times n}$ is a Metzler matrix (resp., a nonnegative matrix). Then $A$ is Hurwitz stable (reps., Schur stable) if and only if there exists a positive vector $\zeta \gg 0$ such that $A\zeta \ll 0 $ (reps., $A\zeta \ll \zeta$). 
\end{lemma}

Finally, for a continuous function $\psi:\R \rightarrow \R$ the upper right Dini derivative is defined as
\begin{equation}\label{dini}
D^+\psi(t)= \limsup_{\delta\rightarrow 0^+}\frac{\psi(t+\delta)-\psi(t)}{\delta}.
\end{equation}
Useful properties of the Dini derivatives can be found in \cite{Rouche}, Appendix 1). 


\section{Main Results}

Consider the  time-varying nonlinear system \eqref{TVS}-\eqref{sector2}. 
By continuity of $f$, it is obvious that $f(0)=0$ and therefore the  system \eqref{TVS} has the zero solution $x(t)\equiv 0, t\geq 0$,  for any $f\in K[\delta,\beta]$. Moreover, as a convention, the function $f\in K[\delta,\beta]$ is assumed to satisfy some Lipschitz or differentiability conditions to assure the global existence and uniqueness of the solution of \eqref{TVS}, for any initial condition. In what follows such a function $f$ is called admissible nonlinearity.

\begin{definition}\label{AES}  The zero solution $x(t)\equiv 0$ of the time-varying nonlinear system \eqref{TVS} is said to be absolutely exponentially stable (shortly, AES) if there exist positive numbers $M, \lambda $  such that for any admissible nonlinearity $f\in K[\delta,\beta] $ and any non-zero continuous function $ \varphi\in C([-h,0],\R^n)$ the solution $x(t)= x(t,\varphi) $  of \eqref{TVS} with the initial condition $ x(\theta)= \varphi(\theta), \ \theta\in [-h,0], $
	satisfies
	\begin{equation}
		\label{conditiondef1}
		\|x(t)\|=	\|x(t,\varphi)\|\leq M e^{-\lambda t}\|\varphi\|, \quad \forall t\geq 0.
	\end{equation}   Such a number $\lambda $ is called the exponential decay rate of $x(t)$. 
\end{definition}

\noindent  Define the Meztler matrix function $	\widehat A(t) =(\widehat a_{k,ij}(t)),\ t\geq 0,\ k\in \underline N,$ by setting 
\begin{equation}\label{metzler}
  \widehat a_{i,i}(t)\!=\! a_{ii}(t),\ \widehat a_{ij}(t)\!=\! |a_{ij}(t)|, j\not=i,\ i,j\in \underline n.  
\end{equation}
The main contribution of this paper is the following
\begin{theorem}\label{main1} Consider the time-varying system \eqref{TVS} with  admissible nonlinearities $f\in K[\delta,\beta]$ 
	where $ \beta = (\beta_1,\ldots,\beta_n)^{\top}$,  $\delta:=(\delta_1,\ldots,\delta_n)^{\top} \in \R^n$ are given  positive vectors such that $\beta_i\geq \delta_i>0,\ \forall i\in \underline n$. Assume that there exists a nonnegative $n\times n$-matrix $\bar B =( \bar b_{ij})  \geq 0$ such that
	\begin{equation}\label{barB}
		|B(t)|\leq \bar B, \ \forall t\geq 0.
	\end{equation}
	Then the zero solution of the  system \eqref{TVS} is AES if there exist  $n$-dimensional vector $\xi :=(\xi_{1},\xi_{2},\ldots,\xi_{n})^{\top}\gg 0$ and a real number $\alpha >0$  
	such that 
	\begin{equation}\label{cond1}
		\bigg( D_{\delta}	\widehat  A^{\top}(t)+ e^{\alpha h}D_{\beta}\bar B^{\top}\bigg)\xi \leq-\alpha \xi, \forall t \geq 0,
	\end{equation}
where $D_{\delta}$ and $D_{\beta}$ are diagonal matrices defined, respectively, as $	D_{\beta}= \text{\rm diag}(\delta_1, \delta_2, \ldots, \delta_n),
	D_{\beta}= \text{\rm diag}(\beta_1, \beta_2, \ldots, \beta_n)$. Moreover, in this case, the exponential decay rate is $\alpha$.
\end{theorem} 
\textit{Proof.} Let $x(t)$ be the solution of \eqref{TVS}, with a nonlinearity  $f\in K[\delta,\beta]$ and a non-zero initial function $\varphi\in C([-h,0],\R^n). $ Then,  for each $i\in \underline n$,  
\begin{equation}\label{mean}
	\dot x_i(t)\! =\! \sum_{j=1}^na_{ij}(t)f_j(x_j(t))\!+\! \sum_{j=1}^nb_{ij}(t)f_j(x_j(t\!-\!h)), \ t\geq 0.
\end{equation} Assume that \eqref{cond1} holds for some $\alpha>0, \xi\gg 0 $. Then, it implies readily
\begin{equation}\label{aii}
\sum_{i=1}^n \widehat a_{ij}(t)\xi_i < 0,\ \forall t\geq 0, \ \forall j\in \underline n. 
\end{equation}
In order to make use of the comparison principle,  let us define the numbers
\begin{equation}\label{d1d2}
 d_1 := \min_{i\in \underline n}\xi_{i},\ \   \ d_2:= \max_{ i\in \underline n}\xi_i.
\end{equation}
and the continuous nonnegative functions $v_0(t), v_1(t),$  by setting,  
\begin{equation*}
	v_0(t)= \sum_{i=1}^n\xi_{i}|x_i(t)|= \xi^{\top}|x(t)|, \ \text{for}\ t\geq -h,
	\end{equation*}
\begin{equation*}
v_1(t) = \sum_{i=1}^n\xi_{i}|\varphi_{i}(t)| \!+\!e^{\alpha h} \sum_{i=1}^n \sum_{j=1}^n \int_{-h}^0e^{\alpha \theta}\bar b_{ij} \xi_{i}\beta_j|\varphi_j(\theta)|d\theta, 
\end{equation*}
for $ t\in [-h,0]$ and 
\begin{equation*}
v_1(t) = v_0(t)+ e^{-\alpha(t-h)} \sum_{i=1}^n \sum_{j=1}^n \int_{t-h}^te^{\alpha s}\bar b_{ij}\xi_i\beta_j|x_j(s)|ds, 
\end{equation*}
for $t\ge 0$. Define, moreover, the following comparison function 
\begin{equation}\label{v1}
	y_{\alpha}(t) = L_1e^{-\alpha t}\|\varphi\|,  \ \text{for}\  t\geq -h, 
\end{equation}	
where  $L_1$ is an arbitrary positive number satisfying
	\begin{equation}\label{Lrho} 
		 d_2 +he^{\alpha h}\sum_{i=1}^ n{\sum_{j=1}^n \bar{b}_{ij}\;\xi_i\;\beta_j} < L_1.	
	\end{equation}
Then, we have obviously
\begin{equation}\label{ineq}
	d_1\sum_{i=1}^n|x_i(t)|\leq v_0(t)\leq v_1(t),   \ \forall t\geq 0.
\end{equation}
Moreover, it follows immediately from the definition of $v_0,v_1,y_{\alpha}$ and \eqref{d1d2}, \eqref{Lrho}  that
\begin{equation} \label{defk} v_1(t)
<  L_1\|\varphi\| \leq  y_{\alpha}(t), \ \text{for} \ t\in [-h,0]. 
\end{equation}
Our goal is to prove that the above inequality holds true for all $t\geq 0$, that is
\begin{equation}\label{goal_0} 
	v_1(t)\leq y_{\alpha}(t) =  L_1e^{-\alpha t} \|\varphi\|, \ \forall t\geq 0. 
\end{equation}
To this end, taking an arbitrary positive number  $\alpha_{\epsilon}\in (0,\alpha)$, we will prove that 
\begin{equation}\label{goal} 
	v_1(t)\leq y_{\alpha_{\epsilon}}(t) =  L_1e^{-\alpha_{\epsilon} t} \|\varphi\|, \ \forall t\geq 0. 
\end{equation}
Assume to the contrary that
\eqref{goal} does not hold. This implies that the set $T_{0}:=\{t\in (0,+\infty): v_{1}(t)> y_{\alpha_{\epsilon}}(t)\}$ is nonempty. Then, denoting $\bar t_0= \inf \{ t\in T_{0} \}$, we have, by the continuity and \eqref{defk},  that $\bar t_0>0 $ and 
\vspace{-0.2cm}
\begin{equation}\label{equal}
	v_1(t)\leq y_{\alpha_{\epsilon}}(t),\ \forall t\in [-h,\bar t_0),\   v_1(\bar t_0)= y_{\alpha_{\epsilon}}(\bar t_0),
\end{equation}
and  there exist a sequence $t_k\downarrow \bar t_0$ such that
\begin{equation}\label{io}
	v_1(t_k) >y_{\alpha_{\epsilon}}(t_k), k=1,2, \ldots 
\end{equation}
Since $f\in K[\delta,\beta] $, it follows immediately from \eqref{sector2} that, for each $i,j\in \underline n$ and $j\not= i$ we have
\begin{equation}\label{deltabeta}
\begin{split} 
\delta_i|x_i|\leq f_i(x_i)\ \text{sign} x_i = |f_i(x_i)| \ \ \text{and}  \\
f_j(x_j)\ \text{sign} x_i\leq |f_j(x_j) |\leq \beta_j|x_j|. 
\end{split}
\end{equation}
Therefore, by using \eqref{metzler}, \eqref{barB}, \eqref{mean}, \eqref{deltabeta},  we get, for each $ t\in [0,\bar t_0], $
\begin{align*} 
	& D^+v_0(t)= \sum_{i=1}^n \xi_{i}D^+ |x_i(t)| \leq  \sum_{i=1}^n \xi_{i} \text{sign}(x_i(t))\dot x_i(t)  \notag \\ 
	& \leq \sum_{i=1}^n\xi_{i}\bigg(a_{ii}(t)|f_i(x_i(t))|+ \sum_{j\not=i}^n |a_{ij}(t)|\;|f_j(x_j(t))|\bigg) \\
&	+ \sum_{i=1}^n\xi_{i}\sum_{j=1}^n|b_{ij}(t)|\;|f_j(x_j(t-h))|\notag\\
	&=\sum_{j=1}^n \sum_{i=1}^n\bigg( \xi_{i}\widehat a_{ij}(t)|f_j(x_j(t))| + \xi_i|b_{ij}(t)||f_j(x_j(t-h))|\bigg)\\
	&\stackrel{\eqref{aii} }\leq \sum_{j=1}^n \sum_{i=1}^n\bigg( \xi_{i}\widehat a_{ij}(t)\delta_j|x_j(t)| + \xi_i \bar{b}_{ij} \beta_j|x_j(t-h)| \bigg). 
\end{align*}Consequently, by the definition of $v_1$ and the properties of the Dini derivative $D^+$ (see, e.g. \cite{Rouche}, Appendix 1) we can deduce, for each $ t\in [0,\bar t_0], $  \vspace{-0.2cm}
\begin{align*} 
 &D^+v_1( t)\!\leq\!\sum_{j=1}^n \sum_{i=1}^n\bigg( \xi_{i}\widehat a_{ij}(t)\delta_j|x_j(t)| \!+\! \xi_i \bar{b}_{ij} \beta_j|x_j(t-h)| \bigg) \\
&+(-\alpha)\left(v_1(t)-v_0(t)\right) +\\
& +e^{\alpha h} \sum_{i=1}^n  \sum_{j=1}^n \xi_{i} \bar b_{ij}\beta_j|x_j(t)| - \sum_{i=1}^n \sum_{j=1}^n \xi_{i}\bar b_{ij}   \beta_j|x_j(t\!-\!h)| \notag\\
    &=-\alpha v_1(t)+\alpha v_0(t)\\
    &+\sum_{j=1}^n \left[\left(D_{\delta}\widehat{A}^{\top}(t)+e^{\alpha h}D_{\beta}\bar{B}^{\top}\right)\xi\right]_j|x_j(t)|\\
&\stackrel{\eqref{cond1}}\leq -\alpha v_1(t)+\alpha v_0(t)+\sum_{j=1}^n (-\alpha)\xi_j|x_j(t)| =-\alpha v_1(t).
\end{align*}
By virtue of  \eqref{equal}, the last inequality implies that\begin{align*}
D^+v_1(\bar{t}_0)&\leq -\alpha v_1(\bar t_0)= -\alpha y_{\alpha_{\epsilon}}(\bar{t}_0) \notag\\
&<-\alpha_{\epsilon} y_{\alpha_{\epsilon}}(\bar{t}_0)= \frac{d}{dt}y_{\alpha_{\epsilon}}(\bar t_0).
\end{align*}
On the other hand, by definition of  $D^+$  and \eqref{equal}, \eqref{io}, we have
\begin{align*} D^+v_1(\bar t_0)
&\geq  \limsup_{t_k\downarrow \bar t_0}\dfrac{v_1(t_k)- v_1(\bar t_0)}{t_k-\bar t_0}\notag\\
&\geq \limsup_{t_k\downarrow \bar t_0}\dfrac{y_{\alpha_{\epsilon}}(t_k)-y_{\alpha_{\epsilon}}(\bar t_0)}{t_k-\bar t_0}=\dfrac{d}{dt}y_{\alpha_{\epsilon}}(\bar t_0),
\end{align*}
conflicting with the above strict inequality. Thus, \eqref{goal} is proved. Now, letting  $\alpha_{\epsilon} \uparrow \alpha$  in \eqref{goal}, we obtain \eqref{goal_0}, 
 which together with \eqref{ineq} implies that 
\begin{equation}\label{tau_0}
\|x(t)\| \leq \frac{v_0(t)}{d_1} \leq Me^{-\alpha t}\|\varphi\|, \ \forall t \geq 0,
\end{equation}
where $ M:= L_1/d_1  $
Therefore, \eqref{conditiondef1} holds, for any $ \varphi \in  C([-h,0],\R^n)$ and any admissible nonlinearity $f\in K[\delta,\beta]$, completing the proof.

As an immediate consequence of Theorem \ref{main1} we obtain the following criterion of AES for time-varying nondelay systems. 

\begin{corollary}\label{cor-Persid} Consider the time-varying nonlinear system
	\begin{equation}\label{nondelay}
		\dot x =A(t)f(x(t)), \ t\geq 0,	
	\end{equation}
	where $f\in K[\delta,\infty)$. Assume that there exist  a positive vector $\xi \gg 0$ such that \vspace{-0.2cm}
	\begin{equation}\label{cond5}
	\gamma := \max_{j\in \underline n} \sup_{t\geq 0}\sum_{i=1}^n\widehat a_{ij}(t)\xi_i \ < \ 0,
	\end{equation} where $\widehat A(t)$ is defined by \eqref{metzler}. Then the zero solution of \eqref{nondelay} is AES, with the exponential decay rate $\alpha,  $ for any $\alpha \in (0, \frac{-\gamma\delta_0}{d_2})$, where $  \delta_0= \min_{j\in \underline n} \delta_j $ and $d_2$ is defined by \eqref{d1d2}.
\end{corollary} 

In the case when both matrix functions $A(t), B(t)$ in \eqref{TVS} can  be upper bounded  by some time-invariant matrices, Theorem \ref{main1} implies the following verifiable criterion of absolute exponential stability.

\begin{corollary}\label{main2} Consider the time-varying  system \eqref{TVS} with  admissible nonlinearities $f\in K[\delta,\beta]$. Assume that there exist a Metzler matrix $\widehat A $, a nonnegative matrix $\bar B \geq 0$ and a positive vector $\xi :=(\xi_{1},\xi_{2},\ldots,\xi_{n})^{\top}\gg 0$  satisfying
	\begin{equation}\label{barAB}
	\widehat A(t) \leq \widehat A,\ 	|B(t)|\leq \bar B, \ \ \forall t\geq 0,
	\end{equation}
and
	\begin{equation}\label{cond2}
		\big( 	\widehat A D_{\delta} + \bar B D_{\beta} \big)^{\top}\xi \ll 0,\  \forall t \geq 0,
	\end{equation}
	where $D_{\delta}$ and $D_{\beta}$ are diagonal matrices defined as in Theorem \ref{main1}. Then the zero solution of system \eqref{TVS} is AES, with the exponential decay rate $\alpha=\alpha_{\max} >0$, which can be calculated, correspondingly to each $\xi \gg 0$  satisfying \eqref{cond2}, as
	\begin{equation}\label{almax}
\alpha_{\max}= \min_{i\in \underline n}\{\alpha_i: g_i(\alpha_i)=0\},
	\end{equation}
where $ g_i(\cdot)$ are continuous functions defined by
\begin{equation}\label{gi}
	g_i(\alpha):= \sum_{j=1}^n\big(\widehat a_{ji}\delta_i\xi_j + e^{\alpha h}\bar b_{ji}\beta_i \xi_j\big)+ \alpha \xi_i, \ i\in \underline n.
	\end{equation}
	\end{corollary}\textcolor{blue}
\noindent {\it Proof.} 
Obviously, for each $ i\in \underline n, \   g_i(\alpha)$ is continuous and monotonically   increasing to $+\infty$ as $\alpha \rightarrow +\infty$  (because $g_i'(\alpha)> 0, \forall \alpha >0$).  Since  $ g_i(0)<0$, due to \eqref{cond2},   it follows that the equation $g_i(\alpha)=0$ has a unique solution $\alpha_i >0 $ and  the inequality $g_i(\alpha)\leq 0$ is valid for all $ \alpha \in [0,\alpha_i]$ but violated for $\alpha>\alpha_i$.  This implies that \eqref{cond1} is valid for all $\alpha \in (0,\alpha_{\max}] $ but violated afterwards. Therefore, by Theorem \ref{main1} the zero solution of the system \eqref{TVS} is AES, with the 'maximal' (for the given $\xi $) exponential decay rate $\alpha_{\max}$.   
 
\begin{corollary}\label{TIDS}
Consider the time-invariant nonlinear system
\begin{equation}\label{persiddelay}
	 \dot x(t) = Af(x(t)) + Bf(x(t-h)), \ t\geq 0,
\end{equation}
where $A$ is a Metzler matrix, $B\geq 0$ and $f$ is any nonlinear function belonging to the sector $ K[\delta,\beta]$, defined by \eqref{sector2}.  Then the zero solution of \eqref{persiddelay} is AES if  there exists a positive vector $\xi \gg 0$ satisfying
\begin{equation}\label{S_H}
(AD_{\delta}+BD_{\beta})^{\top}\xi \ll 0.
\end{equation} 
\end{corollary}
\vspace{0.2cm}
\begin{remark}\label{persidski}
{\rm By Lemma \ref{metzmatrix}, \eqref{S_H} is equivalent to that the Metzler matrix $AD_{\delta}+BD_{\beta}$ is Hurwitz stable. If, conversely, the time-invariant nonlinear system  \eqref{persiddelay} is AES then, since, obviously $f(x):= D_{\delta} x \in K[\delta,\beta]$, it follows that the linear positive system  $\dot x(t)= AD_{\delta}x+BD_{\delta}x(t-h), \ t\geq 0$ is exponentially stable. This, in turn, is equivalent (see, e.g. \cite{Chelaboina2004}) to the existence of a positive vector $\xi \gg 0$ such that $ D_{\delta}(A+B)^{\top}\xi = (AD_{\delta}+BD_{\delta})^{\top}\xi \ll 0$ which implies readily 
\begin{equation}\label{A_S}
	(A+B)^{\top}\xi \ll 0,
\end{equation}
or, equivalently, the Metzler matrix $A+B$ is Hurwitz stable (by Lemma \ref{metzmatrix}). Thus, in the nondelay case (i.e. $B=0$), Corollary \ref{TIDS} is an extension of the main result of \cite{persidskii2005}. 
}	
\end{remark}
\begin{remark}\label{positive} {\rm Note that, under the assumption of Corollary \ref{TIDS}, the nonlinear system \eqref{persiddelay} is \textit{positive} which means that $x(t)=x(t, \varphi)\geq 0, \forall t\geq 0$ for any nonnegative initial function $\varphi \in C([-h,0],\R^n_+)$ (see, e.g. \cite{zhang_zhao2016}, Lemma 4). Therefore, in view of Corollary \ref{TIDS}, Corollary \ref{main_2} amounts to saying that the zero solution of the time-varying system \eqref{TVS}-\eqref{sector2} is AES if its associate  upper bounding (in the sense \eqref{barAB}) by time-invariant positive nonlinear system
			\begin{equation*}
\dot x= \widehat A f(x)+ \bar {B}f(x(t-h), \ t\geq 0
			\end{equation*}  
is AES, under the same  sector constraints. 
}	
\end{remark}

\begin{remark}\label{alex} {\rm 
It is important to mention that  problems of absolute  asymptotic stability have been considered for  time-invariant delay systems of the form \eqref{persiddelay} in a number of previous works, including those with switchings.  
In particular, it has been established (see, e.g. \cite{alex_mason2014,sun_wang2013}) that,  for any nonlinearity $f\in K(0,\infty)$, the solution $x(t)$ of the positive system  \eqref{persiddelay} satisfies $\|x(t)\|\rightarrow 0 $ as $t\rightarrow +\infty$ if there exists $\xi \gg 0$ such that \eqref{A_S} holds.
 This condition, however, is not sufficient to assure the absolute exponential stability of the zero solution (as asserted in Corollary \ref{TIDS} where, however,  a more restrictive nonlinearities class $K[\delta,\beta]\subset K(0,\infty)$ is assumed). Indeed, let us consider the scalar nondelay positive system
$\dot x=-f_0(x)$,  where $f_0(x)=x^3 $ for $|x|\leq 1 $ and 
$f_0(x)=x$ for $ |x|>1.$  Then, clearly, $f_0\in K(0,\infty)$ and \eqref{A_S} holds for $\xi=1$. However, the zero solution $x=0$, while being globally asymptotically stable equilibrium,  is not exponentially stable (e.g., by Theorem 4.6 in  \cite{khalil}, p. 184). Note that  $f_0$ does not belong to $K[\delta,\beta]$ for any $\delta >0$ and any $\beta\geq 1> \delta$.  In \cite{zhang_zhao2016} several criteria of absolute exponential stability for switched time-invariant systems were obtained, also by using the Lyapunov function method, but the conditions look much more complicated  and not easy to be checked. Thus, even for time-invariant systems, our results are novel, while those for the time-varying case as presented in this paper have not yet been known in the existing literature, to the best of our knowledge.}
\end{remark}
\section{Some extensions of the main results}
The approach developed in the previous section can be  extended to get criteria of AES when the system's equation contains time-varying nonlinearities or multiple discrete delays.  We just formulate the results, omitting of proof, because they are largely similar to that for Theorem \ref{main1}.

First, consider a more general model related to the Persidskii-type system \eqref{TVS} that has the form:
\vspace{-0.2cm}
\begin{equation} \label{system_gen}
	\dot{x}_i(t)=\sum_{i=1}^{n}a_{ij}(t)f_{ij}(x_j(t),t)+\sum_{i=1}^{n}b_{ij}(t)f_{ij}(x_j(t-h),t), 
\end{equation}
for $t\geq 0, i=1,\ldots, n$, where the functions $f_{ij}(\cdot, \cdot): \R\times \R_+\rightarrow \R$ are assumed to satisfy, for all $i,j \in \underline n$,
\begin{equation}
	\label{sector_2}
	f_{ij}(x_j, t)x_j>0, \ \forall x_j\ne 0, \text{ and } f_{ij}(0, t)=0, \forall t\geq 0.
\end{equation} 
Then, similarly to Theorem 3.2.10 of \cite{KaszBhaya}, we have the following extension of Theorem \ref{main1}. 

\begin{theorem}\label{main_2} Consider the time-varying system \eqref{system_gen}-\eqref{sector_2} and assume that there exists the admissible nonlinearity  $f=(f_1, f_2, ..., f_n)\in K[\delta,\beta]$ which is defined as \eqref{sector2} such that, for all $i,j\in \underline n$, the following diagonal dominance type conditions  are satisfied 
\begin{equation*}
		\left|f_{ij}(x_j,t)\right|\leq \left|f_{j}(x_j)\right | \leq \left|f_{jj}(x_j,t)\right|,\ \forall  \ i\ne j, \  \forall x_j\in \R,\ \forall t\geq 0.
	\end{equation*}
	Then the zero solution of the  system \eqref{system_gen}-\eqref{sector_2} is AES if there exist  $n$-dimensional vector $\xi :=(\xi_{1},\xi_{2},\ldots,\xi_{n})^{\top}\gg 0$ and a real number $\alpha >0$ 
	such that \eqref{cond1} holds.
\end{theorem} 


Next, consider time-varying nonlinear system with multiple delays of the form \vspace{-0.2cm}
\begin{equation}\label{multiple}
\dot x= A(t)f(x(t))+ \sum_{l=1}^m B_l(t)f(x(t-h_l)),\ t\geq 0,
\end{equation}

\noindent where $A(\cdot), B_l(\cdot), C(\cdot,\cdot)$ are continuous matrix functions and $f\in K[\delta,\beta]$ is an admissible sector-bounded nonlinearity defined by \eqref{sector2}. It is assumed, without loss of generality, that $0<h_1<h_2<\ldots < h_m=h$. Furthermore,  assume that there exist constant matrices $ \widetilde B_l =(\widetilde b_{l,ij})\in \R^{n\times n}$  such that  
\begin{equation}\label{upper}
| B_l(t)| \leq \widetilde{B}_l, \ \forall t \geq 0,\ \forall l \in \underline m.
\end{equation} Then, the following criterion of AES holds for the system \eqref{multiple}.
\begin{theorem}\label{theorem_multiple} Assume that there exist $n$-dimensional vector $\xi=(\xi_1,\xi_2,\ldots,\xi_n)^{\top}\gg 0$ and a real number $\alpha >0$ such that
	\begin{equation}\label{con-multi}
	\bigg( D_{\delta}	\widehat A^{\top}(t)+ \sum_{l=1}^m e^{\alpha h_l}D_{\beta}\widetilde{B}_l^{\top}\bigg)\xi \leq -\alpha \xi,\  \forall t \geq 0,		
	\end{equation} 
\noindent where the matrix function $\widehat A(t)$ and the diagonal matrices $D_{\delta}, D_{\beta}$ are  defined as in Theorem \ref{main1}. Then the zero solution of the delay nonlinear system \eqref{multiple} with sector nonlinearity $f\in K[\delta,\beta]$ is AES. Moreover, in this case,  the exponential decay rate is $\alpha$. 
\end{theorem} 

\section{Time-Varying Nonlinear Difference Systems with delays}

  Consider time-varying nonlinear difference system with delays of the form
\begin{equation}
\label{discrete_system}
x(k\!+\!1)\!=\! A(k)f(x(k))+B(k)f(x(k-h)), \; k\in \mathbb{Z}_+,
\end{equation}
where $h\in \Z_+$ is a given number, $A(k), B(k) : \Z_+ \rightarrow \R^{n\times n}$  are given matrix functions and $f: \R^n\rightarrow \R^n$ is a nonlinear diagonal function belonging to the bounded sector of the form
\begin{equation}\label{sect_disc}
 K(0,\beta]:= \{f: 0<x_if_i(x_i) \leq \beta_ix_i^2, \forall x_i\not=0, i\in \underline n\},
\end{equation}
where $\beta_i >0, i=1,\ldots, n$ are given positive numbers. Such a nonlinear function $f$ is called admissible sector nonlinearity. It is easy to verify that $f\in K(0,\beta]$ if and only if  
\begin{equation}\label{sector_1}
0<x_if_i(x_i), \	0< |f_i(x_i)|\leq \beta_i|x_i|\  \ \text{for}\ \ x_i\not=0, \  i\in \underline n.
\end{equation}
 Let $S[-h,0]$ be the Banach space of functions $\varphi: \mathbb{Z}_{[-h, 0]}\rightarrow \mathbb{R}^n$ equipped with norm 
$ \|\varphi\|=\max_{k\in \mathbb{Z}_{[-h, 0]}}\|\varphi(k)\| $ and $x(k):=x(k, \varphi), k\in \Z_+$ be the solution of  \eqref{discrete_system} satisfying the initial condition $x(k)= \varphi (k), \ \ k\in \Z_{[-h,0]}.$  
We will say that the system \eqref{discrete_system} is absolutely exponentially stable (AES), with the convergence rate $\lambda \in (0, 1) $,  if for a number $L>0$,
\begin{equation}\label{eqdf1}
\|x(k)\|=\left\| x(k, \varphi) \right\| \leq L \;\lambda^k \; \|\varphi\|, \ \forall k\in \mathbb{Z}_+,
\end{equation} for any $\varphi \in S[-h,0]$ and  any admissible sector nonlinearity $f\in K(0,\beta]$. 
\begin{theorem} 
	\label{main_discrete}
	Assume that there exists a vector $ \xi\in \R^n, \xi \gg 0$ and $\lambda\in (0,1)$ satisfying
	\vspace{-0.2cm}
	\begin{equation}\label{pk}
		\bigg(|A(k)|+\lambda^{-h}|B(k)| \bigg) D_{\beta}\xi\leq \lambda \xi,\; \forall k\in \mathbb{Z}_+, \end{equation}
	where $D_{\beta}$ is the diagonal matrix $D_{\beta}= \text{diag}(\beta_1,\ldots, \beta_n)$. 	Then the  system \eqref{discrete_system} is AES with the convergence rate $\lambda$.
\end{theorem}
{\it Proof.} Let $\xi \gg 0$ and $\lambda \in (0,1)$ satisfy \eqref{pk} and $ f \in K(0,\beta]$ be an arbitrary admissible nonlinearity.  Let $x(\cdot)$ be the solution of \eqref{discrete_system} satisfying the initial condition $x(k)= \varphi (k), \ \ k\in \Z_{[-h,0]}.$  Setting $L_0= (\min_{j\in \underline n}\xi_j)^{-1}$ and  $u_i(k)=L_0\xi_i\lambda^k \|\varphi\|,\  k\in \mathbb{Z}_{[-h, +\infty)}, \ i\in \underline{n},$
then we have immediately
\begin{equation}\label{xkpast1}
	|x_i(k)|\leq u_i(k), \ \  \forall k\in \mathbb{Z}_{[-h, 0]}, \forall i\in \underline n.
\end{equation}
Clearly, to prove the theorem, it suffices to show that 
\begin{equation}
	\label{in_discrete1}
	|x_i(k)|\leq u_i(k), \forall k\in  \mathbb{Z}_+, \forall i \in \underline n
\end{equation}
(because then \eqref{eqdf1} holds with $L=L_0\|\xi\|$). Assume to the contrary that \eqref{in_discrete1} does not hold. Then,  in view of  \eqref{xkpast1}, it follows that there exists $k_0>0$, $i_0\in \underline{n}$ such that
\begin{equation}
	\label{in_discrete}
	|x_i(k)|\leq u_i(k), \forall k\in  \mathbb{Z}_{[-h, k_0)}, \forall i\in \underline{n},
\end{equation}
and 
\begin{equation}
	\label{assume_proof}
	|x_{i_0}(k_0)|> u_{i_0}(k_0).\end{equation}
Then, 
we can deduce 
	\begin{align*}		|x_{i_0}(k_0)| 
	&\stackrel{\eqref{discrete_system},\eqref{sector_1}}\leq\sum_{j=1}^n|a_{i_0j}(k_0-1)| \; |x_j(k_0-1)|\beta_j +\sum_{j=1}^n|b_{i_0j}(k_0-1)| \; |x_j(k_0-1-h)|\beta_j\\		&\stackrel{\eqref{xkpast1}, \eqref{in_discrete}}\leq\sum_{j=1}^n |a_{i_0j}(k_0-1)|\xi_jL_0\lambda^{k_0-1}\|\varphi\|\beta_j +\sum_{j=1}^n\lambda^{-h}|b_{i_0j}(k_0-1)| \xi_jL_0\lambda^{k_0-1}\|\varphi\|\beta_j \stackrel{\eqref{pk}} \leq  u_{i_0}(k_0). 
	\end{align*}
This, however, conflicts with \eqref{assume_proof} and completes the proof. 

Theorem  \ref{main_discrete} can be extended to the case of  several delays as follows.
\begin{theorem}\label{th5}
	Consider time-varying nonlinear difference systems with delays of the form
	\begin{equation}
		\label{multdelay_disc}
		x(k\!+\!1)\!=\! A(k)f(x(k)) \!+\!\sum_{l=1}^mB_l(k)f(x(\!k-h_l)), k\in \mathbb{Z}_+,
	\end{equation}
	and $A(\cdot), B_l(\cdot), l\in\underline{ m} $  are given matrix functions on $\Z_+, 0<h_1< ...<h_m$ are given positive numbers and the  nonlinearities $f$ belong to the bounded sector  $ K(0,\beta]$  defined by  \eqref{sect_disc}.
	Then the system \eqref{multdelay_disc} is AES if there exist vector $ \xi \gg 0$ and $\lambda\in (0,1)$ such that,  
	\begin{equation}\label{pk0}
		\bigg(|A(k)|+ \sum_{l=1}^m\lambda^{-h_l}|B_l(k)|\bigg)  D_{\beta}\xi \leq \lambda \xi,\ \forall k\in \Z_+.
	\end{equation}
\end{theorem}
\vspace{0.2cm}
Similarly to the continuous-time case, it is easy to show that for the system  \eqref{multdelay_disc} to be positive, it is necessary and sufficient that all matrices  $A(k), B_l(k)=(b_{l,ij}(k)), k\in \Z_+, l\in \underline m$ are nonnegative. The following  consequence of Theorem \ref{th5} gives a delay-independent criterion of AES for positive difference systems. The proof is similar to that of Corollary \ref{main2}.
\begin{corollary}
	\label{alex_mason}
The delay positive nonlinear difference  system 
	\begin{equation}
		\label{discrete_system2}
		x(k+1)= Af(x(k)) +\sum_{l=1}^mB_lf(x(k-h_l)), \; k\in \mathbb{Z}_+,
	\end{equation}
with sector-bounded nonlinearity  $ f\in K(0,\beta]$ is AES if there exists a vector $ \xi \gg 0$ satisfying 
	\begin{equation}\label{pk1}
		\big(A+B_1+... + B_m\big) D_{\beta}\xi-\xi\ll 0, 
	\end{equation}
Moreover,   in this case, the maximal convergence rate $\lambda_{\max} \in (0,1) $ can be calculated as $\lambda_{\max}:=\max_{i\in \underline{n}}\lambda_i$ where $\lambda_i\in (0,1) $ is the unique solution of the equation $g_i(\lambda)= \sum_{j=1}^n a_{ij}\beta_j\xi_j+ \sum_{l=1}^{m}\lambda^{-h_l}\sum_{j=1}^n b_{l,ij}\beta_j\xi_j-\lambda\xi_i=0.$  Moreover, the above AES property holds  true for 	 {\it any time-varying nonlinear difference system} of the form \eqref{multdelay_disc}, whenever system's matrix functions $A(\cdot), B_l(\cdot), l\in \underline m$ satisfy
	\begin{equation}\label{bound}
	|A(k)|\leq A,\ |B_l(k)|\leq B_l,\ \forall k\in \Z_+, \ \forall l\in\underline m.	
		\end{equation}
\end{corollary}
\vspace{0.5cm}

Corollary \ref{alex_mason} improves considerably the result of \cite{alex_mason2014} (Theorem 6.2) which only proved, equivalently,  that the time-invariant delay positive system \eqref{discrete_system2}, with $\beta=(1, 1, ...,1)^{\top}$, is absolutely {\it asymptotically stable } if \eqref{pk1} holds.  
\vspace{-0.2cm}

\section{Illustrative examples}

\begin{example}\label{Theorem2} {\rm
	We consider the time-varying system of the form \eqref{TVS}, where$n=2$, $h=1$, $f\in K[\delta, \beta]$, with $\delta=(\frac{1}{3}, \frac{1}{2})^{\top}$, $\beta=(\frac{3}{2}, 2)^{\top}$ and for $t\geq 0$,
	\begin{equation*}
		\begin{split}
			&A(t)= \widehat{A}(t)= \begin{bmatrix}-4t-12&0\\t &-2t-5
			\end{bmatrix}, \\
			&B(t)=\begin{bmatrix}\frac{1}{3} \sin t&\frac{1}{8} \cos t\\\frac{1}{3} e^{-t}\cos t&\frac{1}{8} e^{-t}\sin t
			\end{bmatrix}.
		\end{split}
	\end{equation*}
	Clearly, 
		$ |B(t)|\leq \bar{B}= \begin{bmatrix}\frac{1}{3} &\frac{1}{8} \\\frac{1}{3}  &\frac{1}{8} 
		\end{bmatrix}, \forall t\geq 0. $
Then,  taking $\alpha=1$,  $\xi=\begin{bmatrix} 1&1\\
	\end{bmatrix}^{\top}$, we can check that, for all $t\geq 0$, 
	\begin{equation*}
		\left[D_{\delta}\widehat{A}^{\top}(t)+e^{\alpha h}D_{\beta}\bar{B}^{\top}\right]\xi=\begin{bmatrix}-t-4+e\\-t-\frac{5}{2}+\frac{e}{2}
		\end{bmatrix}\ll \begin{bmatrix}-1\\-1
		\end{bmatrix}=-\alpha \xi.
	\end{equation*}
	Then, by Theorem \ref{main1}, 
we conclude that the zero solution the system under consideration is AES, with exponential decay rate $\alpha=1$. 
Note that the above matrix function $A(t)$, $t\geq 0$ can not be upper bounded by any constant matrix, so that  the result of \cite{sun_wang2013, alex_mason2014}  can not be applied in this case.}
\end{example}
\begin{example} {\rm
 Consider the time-varying nonlinear difference system  \eqref{multdelay_disc} where $m=1, h=1, f\in K (0, \beta]$ with $\beta=(\frac{1}{8}; \frac{1}{14})^{\top}$ and, for all  $k\in \mathbb{Z}_+$,
	\begin{equation*}
	A(k)=\begin{bmatrix}-\sin k&2e^{-3k}\\ 3\cos k&-\sin k
			\end{bmatrix}, \ \  B_1(k)=\begin{bmatrix}\frac{1}{2}e^{-k} &\frac{1}{3}\sin k \\ \frac{1}{2} e^{-2k}&\frac{1}{4} \cos k
			\end{bmatrix}, 
	\end{equation*}
	It is easy to see that, for all  $k\in \mathbb{Z}_+$,
	\begin{equation*}
		|A(k)|\leq A=\begin{bmatrix}1&2\\ 3&1
			\end{bmatrix},\ \
		|B_1(k)|\leq B_1=\begin{bmatrix}\frac{1}{2}& \frac{1}{3} \\ \frac{1}{2}&\frac{1}{4}
			\end{bmatrix}.
	\end{equation*}
	Taking $\xi=\begin{bmatrix} 1&1\\ \end{bmatrix}^{\top}$, 
	it can be verified that
	\begin{equation}
		\begin{split}\bigg(A+B_1\bigg) D_{\beta}\xi=\begin{bmatrix}
				\frac{3}{16}+\frac{7}{42}\\ 
				\frac{7}{16}+\frac{5}{56}\\ 
			\end{bmatrix}\leq \begin{bmatrix}
				1\\ 
				1\\ 
			\end{bmatrix}= \xi.
		\end{split}
	\end{equation}
	Therefore, by  Corollary \ref{alex_mason}, the system \eqref{multdelay_disc} is AES and, moreover, the 'maximal' convergence rate is $\lambda_{\max}=0.5840213813$.  }
\end{example} 

\vspace{-0.3cm}

\section{Conclusion}
We have presented a number of verifiable sufficient conditions of absolute exponential stability for different classes of delay time-varying nonlinear systems with sector-bounded nonlinearity. Differently from the traditional approach which is based on using  Lyapunov-Krasovskii functionals, our analysis makes use of comparison principle and the stability characterization of positive upper bounding systems. The results have been obtained for both the continuous-time and discrete-time cases.  When applied to the time-invariant systems, the obtained results have been shown to cover and improve the stability criteria in existing literature. There are several possibilities for developing the work described here. In particular, it would be of interest to investigate  whether our comparison analysis can be adapted to deal with systems with  time-varying delays or to address the absolute exponential stability of time-varying Luri'e  systems. The application of the results obtained in this paper for studying absolute exponential stability of time-varying  nonlinear switched systems would be also an interesting and promising topic which is currently under our consideration. 
\section*{Acknowledgment} This work is supported partly by Vietnam Academy of Science and Technology, via the research project DLTE00.01/22-23.

\vspace{-0.5cm} 



\begin{thebibliography}{00}

\bibitem{lurie1944}A.I. Luri'e,  V.N. Postnikov,  ''On stability theory for controllable systems,'' \emph{ Prikl. Mat. Mekh.},
 vol. 8, no. 3, pp. 246-248, 1944.
 
\bibitem{liberzon_automat2008}
 M.R. Liberzon, "Essays on the absolute stability theory,"  \emph{Autom. Remote Control},  vol. 67, no. 10,
 pp. 1610-1644, 2006.

\bibitem{liao_book} 
X. Liao, P. Yu, \emph{  Absolute Stability of Nonlinear Control Systems}. Springer Science \& Business Media, London, 2008.
 
 \bibitem{KaszBhaya}
 E. Kaszkurewicz, A. Bhaya,\emph{ Matrix Diagonal Stability in Systems and Computation}. Birkhauser, London, 2000.	

 \bibitem{barbashin}
E. Barbashin, "On construction of Lyapunov functions for nonlinear systems,"  in \emph{Proc. 1st IFAC World Congr.,} 1961, pp. 742-751.

\bibitem{persidskii69}
S.K. Persidskii, "Problem of absolute stability," \emph{ Autom. Remote
Control}, vol.12, pp. 1889-1895, 1969.

\bibitem{K_B-SIAM1993}
E. Kaszkurewicz, A. Bhaya, "Robust stability and diagonal Lyapunov functions," \emph{ SIAM J.  Matrix Anal.
Appl.}  vol. 14, no.2,  pp. 508-520, 1993.



\bibitem{Oliveira2002}
M.C. De Oliveira, J.C. Geromel, L. Hsu, "A new absolute stability test for systems with state-dependent perturbations," \emph{ Inter. J. Robust Nonlinear Control}, vol. 14, no. 12, pp. 1209-1226, 2002.

\bibitem{persidskii2005}
S.K. Persidskii, "On the exponential stability of some nonlinear systems," \emph{ Ukrainian Math. J.}, vol.57, pp. 157-164, 2005.




\bibitem{Efimov2021}
D. Efimov, A. Aleksandrov, "On analysis of Persidskii systems and their implementations using LMIs," \emph{ Automatica}, vol. 134, 109905, 2021.





\bibitem{sun_wang2013}
Y. Sun, L. Wang,  "On stability of a class of switched nonlinear systems," \emph{Automatica}, vol. 49, no. 1,  pp. 305-307,  Jan. 2013. 


\bibitem{alex_mason2014}
A. Aleksandrov, O. Mason, " Absolute stability and Lyapunov-Krasovskii functionals for switched nonlinear systems with time-delay," \emph{ J. Franklin Inst.}, vol. 351, no. 8,  pp. 4381-4394, Aug. 2014. 



\bibitem{zhang_zhao2016}
J. Zhang, X. Zhao, J. Huang,  "Absolute exponential stability of switched nonlinear time-delay systems," \emph{ J. Franklin Inst. }, vol. 353, no. 6, pp. 1249-1267, April 2016. 







\bibitem{alex2021}
A. Aleksandrov,  "On the existence of a common Lyapunov function for a family of nonlinear positive systems," \emph{  Syst. Control Lett.}, vol. 147, 104832, Jan. 2021.

\bibitem{zhang2012}
X. Zhao, L. Zhang, P. Shi, M. Liu, "Stability of switched positive linear systems with average dwell time switching,"  \emph{ Automatica},  vol. 48, pp. 1132-1137, 2012.

\bibitem{Liu2018} X. Liu, Q. Zhao, S. Zhong,  "Stability analysis of a class of switched nonlinear systems with delays: a trajectory-based comparison method," \emph{Automatica }, vol. 91, pp. 36-42, May 2018. 

\bibitem{Tian_Sun}
Y. Tian, Y. Sun,  "Exponential stability of switched nonlinear time-varying systems with mixed delays: Comparison principle," \emph{ J. Franklin Inst. }, vol. 357, no.11, pp. 6918-6931, July 2020. 

\bibitem{son_ngoc_2022}
S. Nguyen Khoa, V.N. Le,  "Exponential stability analysis for a class of switched nonlinear time-varying functional differential systems", \emph{ Nonlinear Analysis: Hybrid Systems}, vol.44, 101177, 2022. 







\bibitem{Chelaboina2004}	W.M.  Haddad, V. Chellaboina, "Stability theory for nonnegative and compartmental dynamical systems with time delay," \emph{Syst. Control Lett.}, vol. 51, pp. 355-361, 2004.

\bibitem{Ngoc_IEEE}
P. H. A. Ngoc, "Stability of positive differential systems with delay," \emph{ IEEE Trans. Autom. Control, }, vol. 58, no. 1,  pp. 203-209, Jan.  2013. 

\bibitem{blanc_valcher}
F. Blanchini, P. Colaneri, E. Valcher,  "Switched positive linear systems," \emph{Foundations and Trends
in Systems and Control},  vol. 2, pp.101-273, 2015.

\bibitem{Berman}
A. Berman, R. J. Plemmons, \emph{ Nonnegative Matrices in Mathematical Sciences}. Academic Press, New York, 1979.

\bibitem{Rouche} N. Rouche, P. Habets, M. Laloy, \emph{Stability Theory by Lyapunov Direct Method}. Springer Verlag, Berlin, 1977.

\bibitem{khalil}
H. Khalil,  \emph{ Nonlinear Systems.} Second Ed., Prentice-Hall, Inc., Englewood-Cliffs, NJ, 1996.










\end{thebibliography}
\end{document}